
\input amstex
\input amsppt.sty

\magnification=\magstep1

\define \OO{\Cal O}

\define \C{\Bbb C}

\define \DD{\Bbb D}
\define \FF{\Bbb F}
\define \GG{\Bbb G}
\define \V{\text{V}}
\define \reg{^{\text{reg }}}
\define \Xreg{X\reg}
\define \Sing{\text {Sing }}
\define \IX{\text{I}_X}

\define \Ker{\text{Ker }}
\define \FX{\FF_X}

\define \iso{\cong}
\define \sub{\subseteq}

\define \cn{(\C^n,0)}
\define \w{\omega}

\define \TT{{\roman T}}
\define \T2{\tau}

\define\Def#1{\subheading{Definition}({\it #1})}
\define \proof{\demo{Proof}}
\define \prop #1{\subheading{Proposition}({\it #1})\par}
\define \theo #1{\subheading{Theorem}({\it #1})\par}
\define \lemma #1{\subheading{Lemma}({\it #1})\par}
\define \abs{\medpagebreak\flushpar}

\topmatter
\title Conormal Differential Forms of an Analytic Germ
\endtitle

\abstract A differential form vanishing on the tangent space at
smooth points of a reduced embedded analytic germ is called
conormal. For proving that a conormal one--form of a hypersurface
vanishes at its singularities we state a Bertini--type theorem.
\endabstract

\author Robert Gassler
\endauthor

\address
Department of Mathematics, 567 Lake, Northeastern University, Boston, MA 02115,
U.S.A., e-mail: gassler\@neu.edu
\endaddress
\endtopmatter

\heading{0. Introduction}
\endheading
\abs The purpose of this paper is to study properties of
differential forms which vanish on the tangent space of an
analytic germ at its regular points. We will call them {\it
conormal} to the germ. They form a differential ideal the
properties of which will be studied in section 2. In section 4 we
state our main result saying that conormal one--forms of a
hypersurface vanish at its singularities. For the proof we need
two ingredients: a Bertini--type theorem which we prove in section
3, and a result about tangential vector fields which are
isomorphic to conormal $(n-1)$--forms. The last section gives some
examples of our forms.

This work was finished during my stay in Valladolid, Spain,
financed by an ERASMUS--scholarship. I want to thank the members
of the Department of Algebra, Geometry, and Topology for their
hospitality, and of course Herwig Hauser for supervising this work
and spending his time on dicussing problems with me.

\heading{1. Notation}
\endheading

\abs The germs of holomorphic functions on $\C^n$ at 0 will be
denoted by $\OO,$ holomorphic differential forms by $\Omega,$
holomorphic vector fields by $\DD.$ All our analytic germs will be
embedded and reduced. We will identify them with their
well--chosen  representatives. In the following $X\sub\cn$ will be
an analytic germ. The germ of regular points of $X$ is denoted by
$\Xreg,$ the singular locus by $\Sing X,$ and the Lie--Algebra of
tangential vector fields by $\DD_X.$

\abs The Zariski--tangent space of $X$ at 0 is denoted by
$\TT_0X.$  In the following we will use the identification
$\TT_0\C^n\iso\C^n$ induced by the canonical chart on $\C^n$
without mention.

The {\it tangent plane space $\T2 X$ of $X$} is the closure of the
tangent bundle $\TT\Xreg\sub \cn\times \C^n.$  If $X$ is of pure
dimension $r$ the closure of
$\{(p,\TT_pX),p\in\Xreg\}\sub\cn\times\GG^{n-1,r-1}$ will also be
called the {\it tangent plane space $\T2^*X$ of $X$}. The fibers
of these bundles over a point $p$ will be denoted by subindex $p.$
The tangent plane space is compatible with the inclusion of
analytic germs: for analytic germs $X\sub Y\sub \cn$ we have $\T2
X\sub \T2 Y,$ as proved by Whitney [Wh, p. 548], as an application
of the existence of regular stratifications.

\heading 2. Conormal Forms of Analytic Germs \endheading

\Def{Conormal differential forms}
A differential form $\w\in\Omega$ is called {\it conormal} to $X,$
if it vanishes on $\TT\Xreg.$ This means each $i$--homogeneous
part of $\w$ vanishes on $(\TT_pX)^i$ for $p\in\Xreg.$

The set of all conormal differential forms of $X$ is denoted by
$\FX,$ the set of $k$--homogeneous conormal forms by $\FX^k.$

\prop{Properties of conormal forms}
\item{(i)} $\FX$ is a differential ideal in $\Omega.$
\item{(ii)} Let $X=\V(f_1,\dots,f_m)$ be a complete intersection.
A differential form $\w \in \Omega$ is conormal to $X$ iff
$\w\wedge df_1\wedge \dots \wedge df_m$ vanishes on $X.$
\item{(iii)}The ideal of $X$ equals $\FX^0.$
\item{(iv)} A differential form is conormal to $X$ iff it
annihilates the tangent plane space $\T2 X.$
\item{(v)} An analytic germ $Y\sub \cn$ is contained in $X$ iff
$\FX\sub \FF_Y.$
\item{(vi)} Let $X_1,\dots,X_m$ be the irreducible components of
$X.$ Then $\FX=\FF_{ X_1}\cap\dots\cap\FF_{X_m}.$
\item{(vii)} The sequence
$$0\to \FX^0\overset d \to\to \FX^1 \overset d \to\to\dots
\overset d \to \to \FX^n \to 0$$
is exact at $\FX^0$ and $\FX^1.$
\demo{Proof}(i) The pull--back of differential forms under the
embedding $i:\Xreg\hookrightarrow\cn$ is a homomorphism of
differential algebras. $\FX$ is its kernel, hence a differential--ideal.

\abs (ii) We must prove the assertion only at regular points,
because $df_1\wedge\dots\wedge df_m$ vanishes on $\Sing X.$
Furthermore, $\w$ can be assumed to be $r$--homogeneous, because
of the linearity of the exterior product. For $r=0$ the assertion
is obvious, since $df_1\wedge\dots\wedge df_m$ does not vanish on
$X.$ For $r>0$ and $p\in \Xreg$ we know that the linear forms
$df_i(p)$ are linearly independent. Hence, we can extend them to a
basis $(l_i)$ of $\TT^*_p X.$  Now $\w(p)$ annihilates
$\TT_pX=\bigcap_1^m\Ker l_i$ iff it is an element of the ideal
generated by $l_1,\dots,l_m$ in the $\C$--algebra
$\bigwedge\TT^*_p\C^n.$ This is equivalent to $\w(p)\wedge
l_1\wedge\dots\wedge l_m=(\w\wedge df_1\wedge \dots \wedge
df_m)(p)=0.$

\abs (iii) This is clear by the definition of $\FX^0.$

\abs(iv) Follows from the definition of $\T2 X.$

\abs(v) For $Y\sub X$ we have $\T2 Y\sub\T2 X$. With (iv) we get
$\FX\sub\FF_Y.$ If $\FX\sub\FF_Y$ we have with (iii)
$$\IX=\FX^0\sub\FF_Y^0=\text{I}_Y.$$

\abs(vi) The inclusion $\sub$ follows from (v). On the other hand
every $\w\in\FF_{X_1} \cap\dots\cap\FF_{X_m}$ vanishes on the
tangent space of $X$ at a regular point, as it is also a regular
point in one irreducible component.

\abs(vii) The only thing to show is that for every $\w\in\FX^1$
with $d\w=0$ there is a $g\in\FX^0=\IX$ with $dg=\w.$ The Lemma of
Poincar\'e gives us the existence of a function $f\in\OO$ with
$\w=df.$ We now show $g\:\!\!= f-f(0)\in\IX.$  For each point
$p\in\Xreg$ we can find a continuous, piecewise differentiable
path $\gamma\:[0,1]\to X$ connecting 0 and $p,$ which intersects
$\Sing X$ in a finite number of points, as $\Xreg$ is a manifold,
dense in $X,$ and has only a finite number of connected
components. We now have  $$g(p)=f(p)-f(0)=\int_\gamma
df=\int_0^1\w(\gamma(t)).\dot\gamma(t)=0,$$ because
$\dot\gamma(t)\in\TT_{\gamma(t)}X\sub\Ker \w(\gamma(t))$ for all
but a finite number of times.
\enddemo

\heading 3. Intersections of Hyperplanes and Hypersurfaces
\endheading

\abs
In this section we will prove a special type of Bertini's Theorem,
which could also be derived from the results of Flenner [Fl] and
Jouanalou [J]. Here, an idea of Teissier [T, Lemma 1.1] will be
further developed to give a very geometric proof of the theorem
for our special case. We will also consider non--reduced analytic
germs in this section.

\abs The following proposition collects results proved by
Teissier [T, Lemma1.1] and Whitney [Wh, Lemma 6.1]. For the
Grassmannian $\GG^{n-1,n-2}$ of hyperplanes in $\C^n$ we will
write $\GG.$ A hyperplane will be called {\it transversal} to $X$
if it is transversal to $\Xreg.$

\prop {Transversal Hyperplanes}
Let $X\sub\cn$ be of pure dimension greater than 0. The
intersection of $X$ and a transversal hyperplane has dimension
$\dim X-1.$ If $X$ is smooth, so is $H\cap X.$ The non--empty  set
$W=\{H\in \GG: H\text{ is transversal to every element of
}\T2^*_0X\}$ is dense in $\GG.$ Each of its elements is
transversal to  $X$.

\theo {Theorem of Bertini for analytic germs}
Let $X$ be a reduced hypersurface in $\cn,n\geq 3.$  For each
hyperplane $H\in\GG$ transversal to $X$ and $\Sing X,$ $H\cap X$
is reduced and $$\Sing(X\cap H)=H\cap\Sing X.$$  By the above
proposition the set of these hyperplanes is dense in $\GG.$

\demo{Proof} For every $H\in \GG$ transversal to $X$ we have
$H\cap\Xreg\sub(H\cap X)\reg$ by the above proposition. On the
other hand, we know that $H\cap\Sing X\sub\Sing(X\cap H).$

If the hypersurface $H\cap X\sub H$ were not reduced, at least one
of its irreducible components would be contained in its singular
locus. Hence, $\Sing X\cap H$ would be one--codimensional in $H.$
But, it has codimension at least two  which is obvious if $\Sing
X$ equals 0, and otherwise follows from the above proposition.
\enddemo

\heading 4. Conormal 1--Forms of a Hypersurface \endheading

We now study the behavior of 1--forms conormal to a hypersurface
at its singularities. In the following $X\sub\cn$ is a
hypersurface, whose ideal is generated by $f\in\OO.$

\prop {Conormal $(n-1)$--Forms of a Hypersurface}  A conormal $(n-
1)$--form of $X$ not vanishing at 0 induces a decomposition $X\iso
X'\times(\C,0),$ where $X'\sub(\C^{n-1},0)$ is an analytic germ.

\demo{Proof} The $\OO$--linear map $D\:\Omega^{n-1}\to\DD,\w
\mapsto(f\mapsto\phi(\w\wedge df)), $ with $\phi$ the isomorphism
mapping $dx_1\wedge\dots\wedge dx_n$ to 1, is a bijection and maps
$\FX^{n-1}$ to $\DD_X,$ as $(D(\w).f)|_X=\phi(\w\wedge df|_X)=0$
for $\w$ conormal to $X.$ Applying the Theorem of Rossi to the
image under $D$ of a form conormal to $X$ and not vanishing at 0
gives us the decomposition
[R, Theorem 3.2].
\enddemo

\subheading{Remark}
Note that we used the natural orientation and metric on $\C^n$ to
define the bijection $\phi\:\Omega^{n-1}\to \OO.$

\theo{Conormal 1--forms of a hypersurface} Each conormal 1--form
of $X$ vanishes on $\Sing X.$

\proof

\abs -- $dim\, X=1$: A 1--form not vanishing in 0 induces the
decomposition $X\iso(\C,0)\times X', X'$ is analytic, as proved in
the above proposition. It is obvious that such a decomposition
only exists if $X$ is regular.

\abs -- {\it Induction on $dim\, X$}: Suppose  the theorem to be
proved for hypersurfaces $X\sub(\C^{n-1},0) $ and look at a
singular point which can be considered to be 0 without loss of
generality. By the Theorem of Bertini there exist two hyperplanes
$H_1,H_2\sub\cn$ transversal to each other such that $H_i\cap X$
is reduced and  $\Sing (H_i\cap X)= H_i\cap \Sing X.$ Consider now
a conormal 1--form $\w$ of $X.$ By (v) of the  proposition of
properties of conormal forms, $\w|_{\TT H_i}$ is also conormal to
the hypersurface $H_i\cap X\sub H_i.$ For this case our theorem is
supposed to be proved, so $\w$ vanishes on $\TT_0H_i.$ These two
tangential spaces span $\C^n$ and we get $\w(0)=0.$
\enddemo
\heading 5. Examples
\endheading

\abs
In this section we will use the following result proved by Vetter
and Lebelt [V],[L]. Furthermore, $\FX^k$ will be called {\it
trivial} if it equals the $\OO$--module of k--homogeneous forms of
the differential ideal generated by $\IX$ in $\Omega.$

\theo{Trivial Conormal Forms}
Let $X\sub\cn$ be a complete intersection. For $k\leq \dim X$ the
following are equivalent:
\item{(i)} $X$ is regular in codimension $k$
\item{(ii)} $\FX^k$ is trivial.

\abs 1. $X=\V(x_1,\dots,x_m)\sub\cn.$ The differential ideal
$\FX$ of $\Omega$ is generated by $(x_1,\dots, x_m).$ This can be
proved by an easy computation using the characterisation of
conormal forms, but it also follows directly from the above
theorem.

\abs 2. $X=\V(x^3-yz).$  $\FX^1$ is
trivial, but there is a non--trivial conormal 2--form $\w=x
dy\wedge dz+ 3zdx\wedge dy,$ which generates, together with the
trivial forms, $\FX^2$ as an $\OO$--module.

\abs 3. $X=\V(z^2-xy^2),$ the {\it Whitney--umbrella}. $\FX$ is
generated by $z^2-xy^2,yzdx+2xzdy-2xydz,$ and $ydx\wedge dz-
zdx\wedge dy.$

\abs 4. $X=\V(xz-yt).$ Both $\FX^1$ and $\FX^2$ are trivial. A
non--trivial 3--form is $\w=(xdy- y dx)\wedge dz \wedge dt.$

\enddocument